\documentstyle[12pt]{article}

\begin{document}

\date{}

\title{\Large\bf Functional Integration and the Kontsevich Integral}

\author{Louis H. Kauffman \\
  Department of Mathematics, Statistics and Computer Science \\
  University of Illinois at Chicago \\
  851 South Morgan Street\\
  Chicago, IL, 60607-7045}

 \maketitle
  
 \thispagestyle{empty}
 
 \subsection*{\centering Abstract}

{\em
 This paper is an exposition of the relationship between Witten's functional integral and Vassiliev invariants. 
}
 
\section{Introduction}

\noindent
This paper shows how the Kontsevich Integrals, giving Vassiliev 
invariants in knot theory, arise naturally in the perturbative expansion of Witten's functional integral.  The paper is a sequel to \cite{WittKont}. Since the writing of \cite{WittKont} I became aware of the work of Labastida and P$\acute{e}$rez \cite{LP} on this same subject. Their work comes to an identical conclusion, interpreting the Kontsevich integrals in terms of the light-cone gauge and thereby extending the original work of Fr\"ohlich and King \cite{Frohlich and King}.  The purpose of this paper is to give an exposition of these relationships and to introduce  diagrammatic techniques that illuminate the connections. In particular, we use a diagrammatic operator method that is useful both for Vassiliev invariants and for relations of this subject with the quantum gravity formalism of Ashtekar, Smolin and Rovelli \cite{ASR}. An aspect that this paper does not treat is the perturbation expansion via three-space integrals leading to Vassiliev invariants as in \cite{Altschuler-Friedel}. See also \cite{Bott-Taubes}.  Nor do we deal with the combinatorial reformulation of Vassiliev invariants that proceeds from the Kontsevich integrals as in \cite{Cart}. 
\vspace{3mm}

The paper is divided into three sections.   Section 2 discusses Vassiliev invariants and 
invariants of rigid vertex graphs.  The section three on the functional integral introduces  
the basic formalism and shows how the functional integral is related directly to Vassiliev 
invariants. In this section we also show how our formalism works for the loop transform of Ashtekar,Smolin and Rovelli.  
Finally section 4 shows how the Kontsevich integral arises in the perturbative 
expansion of Witten's integral in the axial gauge. One feature of section 4 is a new and 
simplified calculation of the necessary correlation functions by using the complex numbers and the two-dimensional Laplacian.   We show how the Kontsevich integrals are the Feynman integrals for this theory.  
\vspace{3mm}

\noindent
{\bf Acknowledgement.}    It gives the author pleasure to thank Louis Licht, Chris King 
and Jurg Fr\"ohlich for helpful conversations and  to thank the National Science 
Foundation for support of this research under NSF Grant DMS-9205277 and the NSA for 
partial support under grant number MSPF-96G-179.
 \vspace{3mm}

\section{Vassiliev Invariants and Invariants of Rigid Vertex Graphs} 

If  $V(K)$ is a  (Laurent polynomial valued,  or more generally - commutative ring valued)  
invariant of knots,  then it can be naturally extended to an invariant of rigid vertex graphs 
\cite{Kauffman-Graph} by defining the invariant of graphs in terms of the knot invariant 
via an unfolding  of the vertex. That is, we can regard the vertex as a "black box" and 
replace it by any tangle of our choice. Rigid vertex motions of the graph preserve the 
contents of the black box, and hence implicate ambient isotopies of the link obtained by 
replacing the black box by its contents. Invariants of knots and links that
 are evaluated on these replacements are then automatically rigid vertex invariants of the 
corresponding graphs. If we set up a collection of multiple replacements at the vertices
 with standard conventions for the insertions of the tangles, then a summation over all 
possible replacements can lead to a graph invariant with new coefficients corresponding to 
the different replacements.  In this way each invariant of knots and links implicates a large 
collection of graph invariants. See \cite{Kauffman-Graph}, \cite{Kauffman-Vogel}. 
 \vspace{3mm}
 
 The simplest tangle replacements for a 4-valent vertex are the two crossings, positive and 
negative, and the oriented smoothing. Let V(K) be any invariant of knots and links.
 Extend V to the category of rigid vertex embeddings of 4-valent graphs by the formula 
  $$V(K_{*}) = aV(K_{+}) + bV(K_{-}) + cV(K_{0})$$
where $K_{+}$ denotes a knot diagram $K$ with a specific choice of positive crossing, 
$K_{-}$ denotes a diagram identical to the first with the positive crossing replaced by a 
negative crossing and  $K_{*}$ denotes a diagram identical to the first with the positive 
crossing replaced by a graphical node. 
  \vspace{3mm}

This formula means that we define  $V(G)$  for an embedded 4-valent graph  $G$  by 
taking the sum
  
    $$V(G) = \sum_{S} a^{i_{+}(S)}b^{i_{-}(S)}c^{i_{0}(S)}V(S)$$ 
   
\noindent
 with the summation over  all knots and links $S$ obtained from  $G$ 
 by replacing a node of $G$ with either a crossing of positive or negative type, or with  a 
smoothing of the crossing that replaces it by a planar embedding of non-touching segments 
(denoted $0$).  It is not hard to see that if $V(K)$  is an  ambient isotopy invariant of 
knots, then,  this extension is an rigid vertex isotopy invariant of graphs.  In rigid vertex 
isotopy the cyclic order at the vertex is preserved, so that the vertex behaves like a rigid 
disk with
 flexible strings attached to it at specific points.   
\vspace{3mm}

There is a rich class of graph invariants that can be studied in this manner.  The Vassiliev 
Invariants   
\cite{Vassiliev},\cite{Birman and Lin},\cite{Bar-Natan}
 constitute the important special case of these graph invariants where  $a=+1$, $b=-1$ and 
$c=0.$    Thus  $V(G)$  is a Vassiliev invariant if

			$$V(K_{*}) = V(K_{+})  -  V(K_{-}).$$

\noindent
 Call this formula the {\em exchange identity} for the Vassiliev invariant $V.$  See Figure 1
\vspace{3mm}

\begin{figure}[htbp]
\vspace*{140mm}
\includegraphics{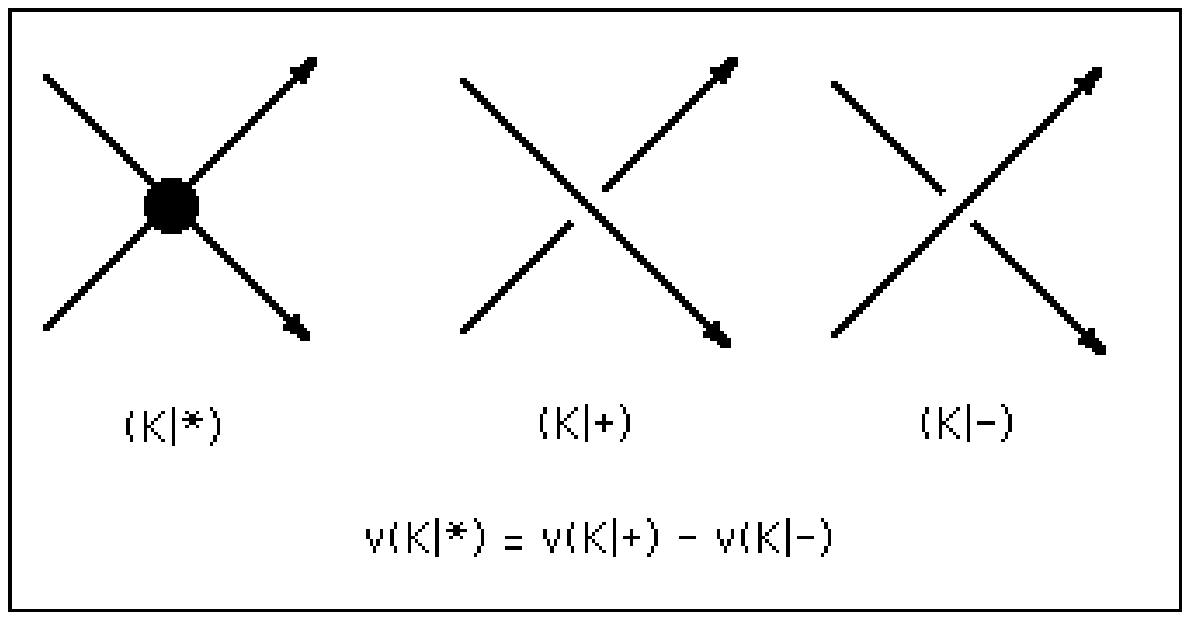}
\vspace*{13pt}
\begin{center}
{\bf Figure 1 --- Exchange Identity for Vassiliev Invariants}
\end{center}
\end{figure}
\vspace{3mm}

 $V$  is said to be of  {\em finite type}  $k$  if  $V(G) = 0$  whenever  $|G| >k$  where 
$|G|$  denotes the number of (4-valent) nodes in the graph $G.$ The notion of finite type is 
of extraordinary significance in studying these invariants. One reason for this is the 
following basic Lemma.
 \vspace{3mm}
 
 \noindent {\bf Lemma.} If a graph $G$ has exactly $k$ nodes, then the value of a 
Vassiliev invariant $v_{k}$ of type $k$ on $G$, $v_{k}(G)$, is independent of the 
embedding of $G$.
 \vspace{3mm}
 
\noindent {\bf Proof.} The different embeddings of $G$ can be represented by link 
diagrams with some of the 4-valent vertices in the diagram corresponding to the nodes of 
$G$. It suffices to show that the value of $v_{k}(G)$ is unchanged under switching of
a crossing.  However, the exchange identity for $v_{k}$ shows that this difference is 
equal to the evaluation of $v_{k}$ on a graph with $k+1$ nodes and hence is equal to 
zero. This completes the proof.//
 \vspace{3mm}
 
 The upshot of this Lemma is that Vassiliev invariants of type $k$ are intimately involved 
with certain abstract evaluations of graphs with $k$ nodes. In fact, there are  restrictions 
(the four-term relations) on these evaluations demanded by the topology  and it follows 
from results of Kontsevich \cite{Bar-Natan} that such abstract evaluations actually 
determine the invariants. The knot  invariants derived from classical Lie algebras are all 
built from Vassiliev invariants of finite type. All this is directly related to Witten's 
functional integral \cite{Witten}.  
 \vspace{3mm}

In the next few figures we illustrate some of these main points.
In Figure 2 we show how one associates a so-called chord diagram to represent the abstract graph associated with an embedded graph. The chord diagram is a circle with arcs connecting those points on the circle that are welded to form the corresponding graph.  In Figure 3 we illustrate how the four-term relation is a consequence of topological invariance. In Figure 4 we show how the four term relation is a consequence of the abstract pattern of the commutator identity for a matrix Lie algebra. This shows that the four term relation is directly related to a categorical generalisation of Lie algebras. Figure 5 illustrates how the weights are assigned to the chord diagrams in the Lie algebra case -  by inserting Lie algebra matrices into the circle and taking a trace of a sum of matrix products.
\vspace{3mm}

\begin{figure}[htbp]
\vspace*{80mm}
\includegraphics{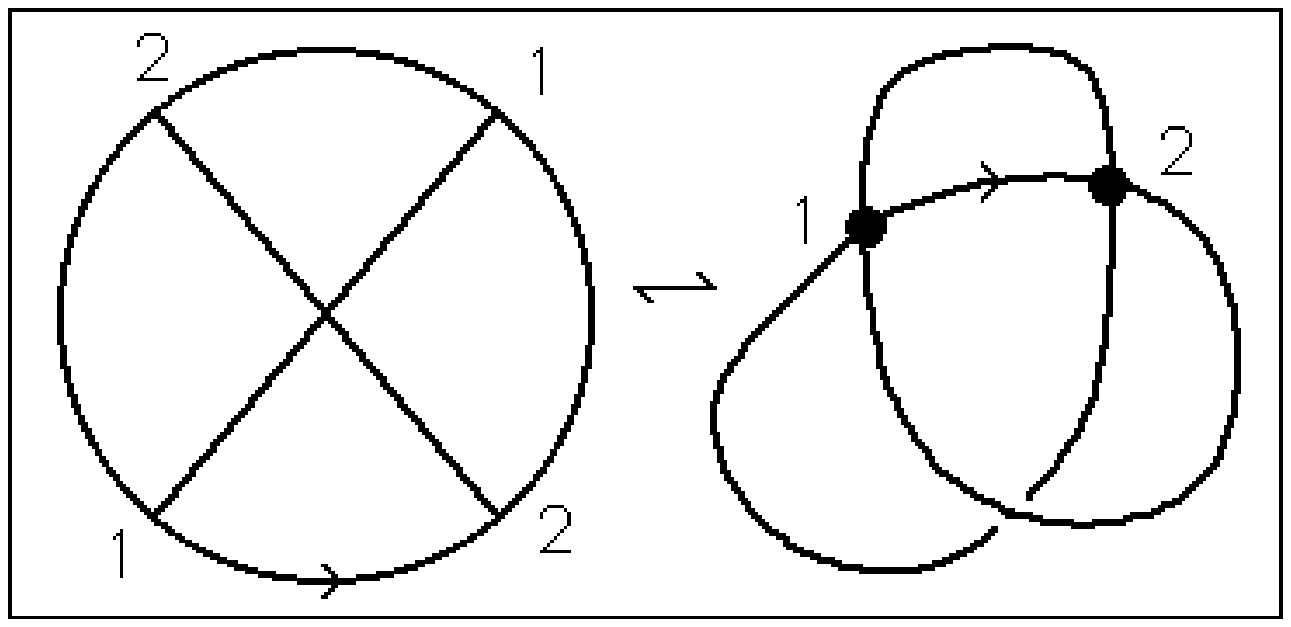}
\vspace*{13pt}
\begin{center}
{\bf Figure 2 ---  Chord Diagrams}
\end{center}
\end{figure}
\vspace{3mm}

\begin{figure}[htbp]
\vspace*{160mm}
\includegraphics{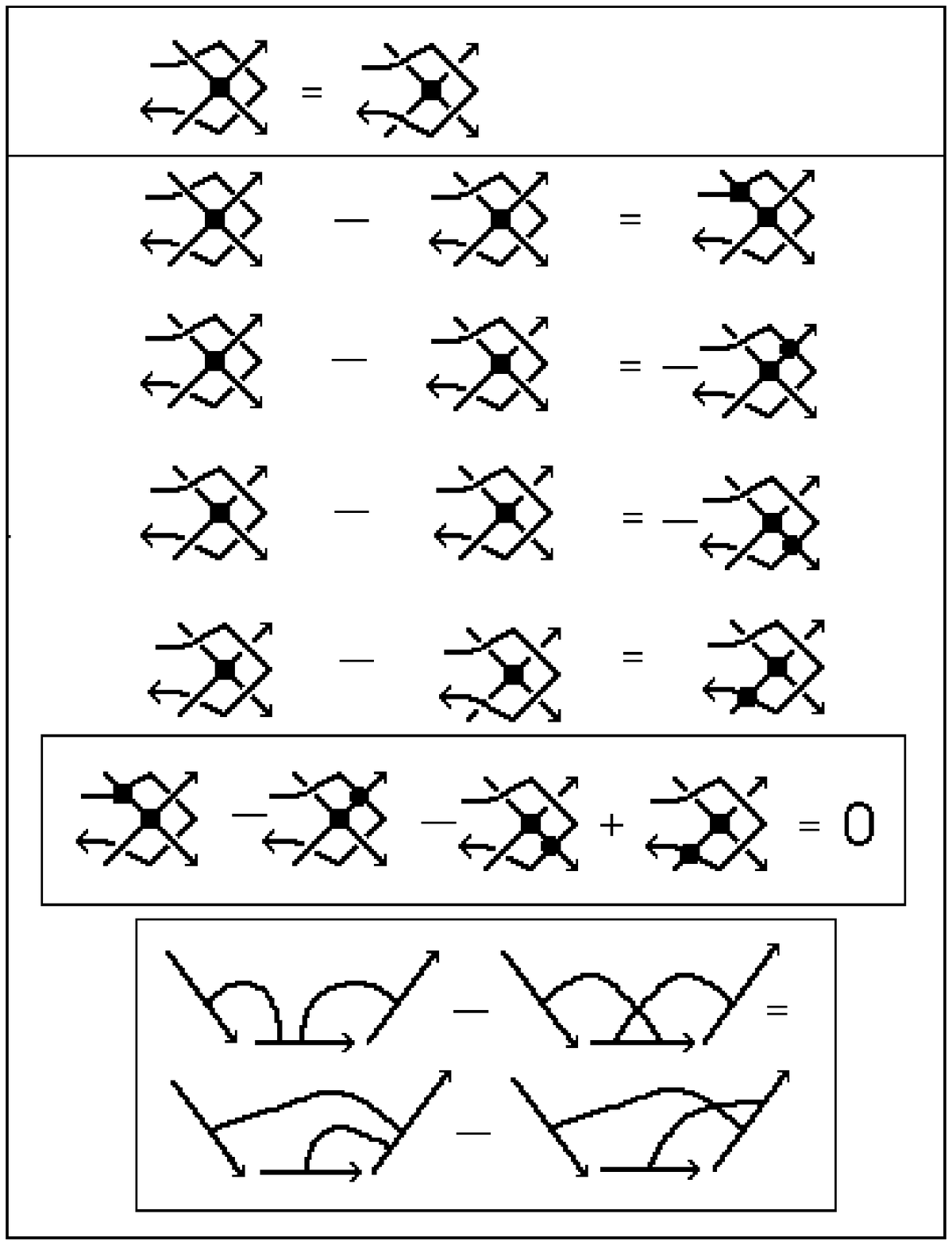}
\vspace*{13pt}
\begin{center}
{\bf Figure 3 --- The Four Term Relation from Topology}
\end{center}
\end{figure}
\vspace{3mm}

\begin{figure}[htbp]
\vspace*{160mm}
\includegraphics{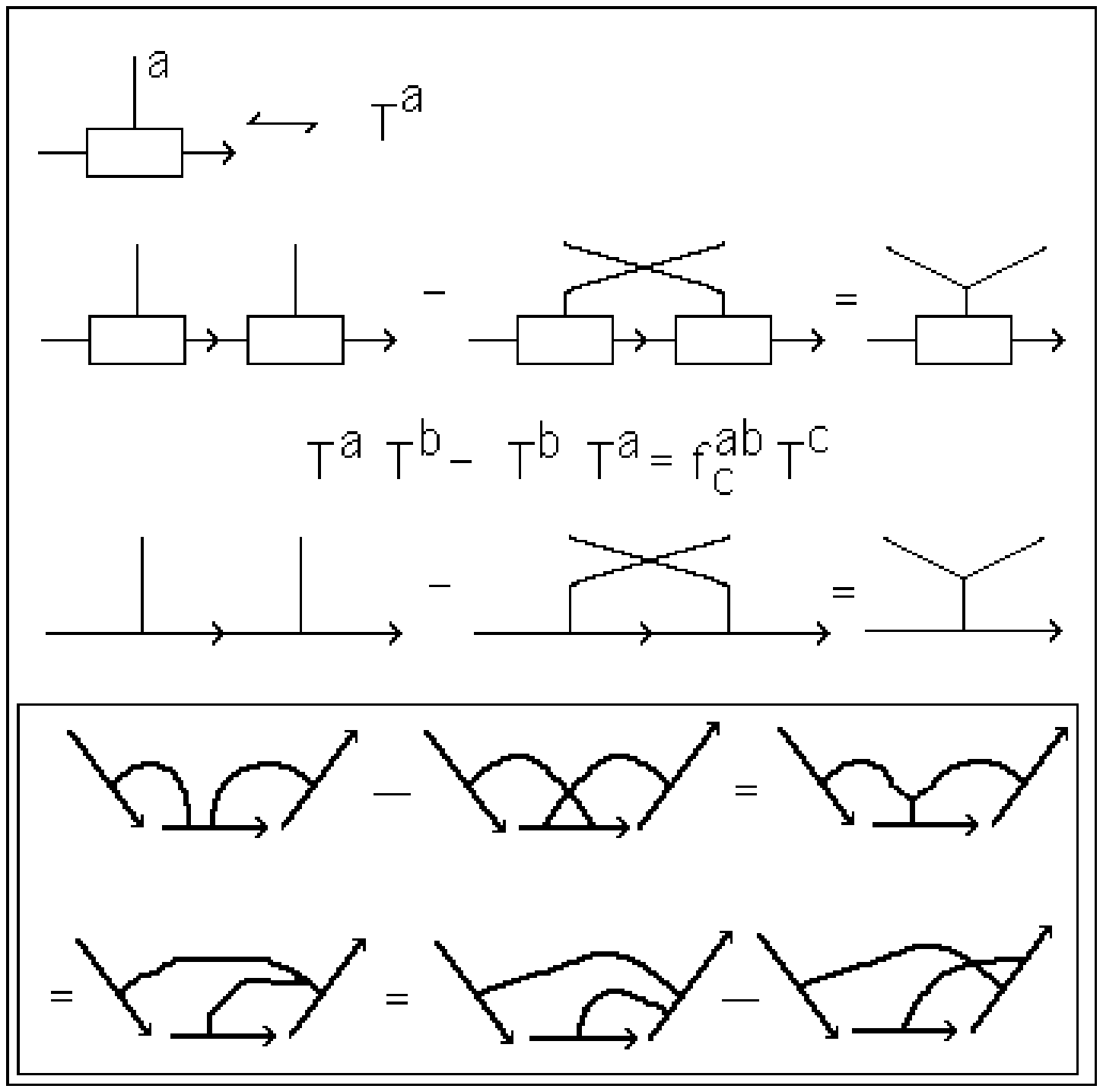}
\vspace*{13pt}
\begin{center}
{\bf Figure 4 --- The Four Term Relation from Categorical Lie Algebra}
\end{center}
\end{figure}
\vspace{3mm}

\begin{figure}[htbp]
\vspace*{60mm}
\includegraphics{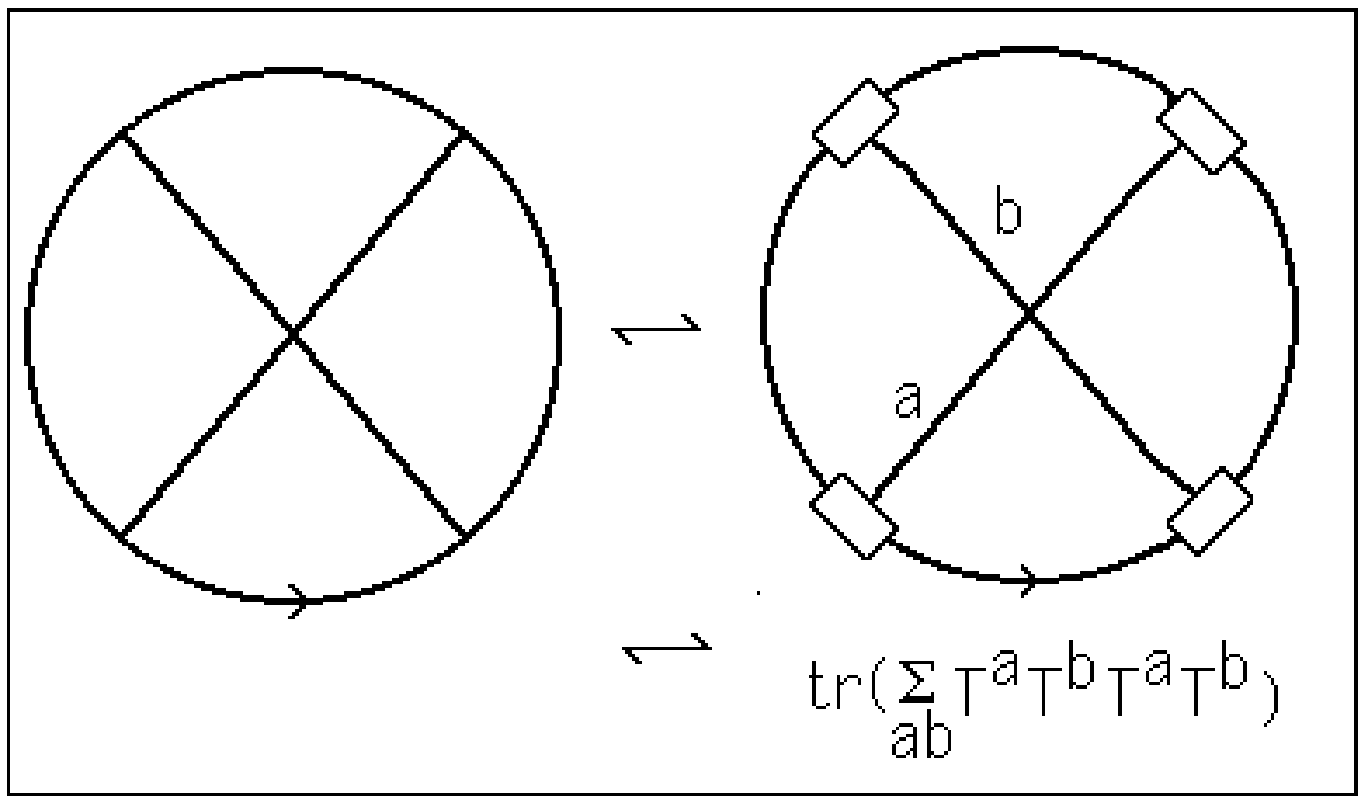}
\vspace*{13pt}
\begin{center}
{\bf Figure 5 --- Calculating Lie Algebra Weights}
\end{center}
\end{figure}
\vspace{3mm}

\section{Vassiliev Invariants and Witten's Functional Integral}

In \cite{Witten}  Edward Witten proposed a formulation of a class of 3-manifold invariants 
as generalized Feynman integrals taking the form  $Z(M)$  where
   
$$Z(M) = \int DAe^{(ik/4\pi)S(M,A)}.$$

\noindent			  
Here  $M$ denotes a 3-manifold without boundary and $A$ is a gauge field  (also called a 
gauge potential or gauge connection)  defined on $M$.  The gauge field is a one-form on a 
trivial $G$-bundle over  $M$ with values in a representation of the  Lie algebra of $G.$ 
The group $G$ corresponding to this Lie algebra is said to be the gauge group. In this 
integral the action   $S(M,A)$  is taken to be the integral over $M$ of the trace of the 
Chern-Simons three-form   $A \wedge dA + (2/3)A \wedge A \wedge A$.  (The product is 
the wedge product of differential forms.)
\vspace{3mm} 

$Z(M)$  integrates over all gauge fields modulo gauge equivalence (See \cite{Atiyah:YM} 
for a discussion of the 
definition and meaning of gauge equivalence.)
\vspace{3mm}   

The formalism  and   internal logic of Witten's integral supports  the existence of a large 
class of topological invariants of 3-manifolds and  associated invariants of knots and links 
in these manifolds.
\vspace{3mm} 

The invariants associated with this integral have been given rigorous 
combinatorial descriptions  \cite{RT},\cite{Turaev-Wenzl},\cite{Kirby-Melvin},\cite{Lickorish}, \cite{Walker},\cite{TL},  
but questions and conjectures arising from the integral formulation are still outstanding. 
(See for example \cite{Atiyah}, \cite{Garoufalidis},\cite{Gompf&Freed}, 
\cite{Jeffrey},\cite{Rozansky}, \cite{Adams}.)
Specific conjectures about this integral take the form of just how it implicates invariants of 
links and 3-manifolds, and how these invariants behave in certain limits of the coupling 
constant $k$ in the integral. Many conjectures of this sort can be verified through the 
combinatorial models. On the other hand, the really outstanding conjecture about the 
integral is that it exists! At the present time there is no measure theory or generalization of 
measure theory that supports it.   Here is a formal structure of great beauty. It is also a 
structure whose consequences can be verified by a remarkable variety of alternative means.  
\vspace{3mm}

We now  look at the formalism of the Witten integral in more detail and see how it 
implicates invariants of knots and links corresponding to each classical Lie algebra.   In 
order to accomplish this task, we need to introduce the Wilson loop.  The Wilson loop is an 
exponentiated version of integrating the gauge field along a loop  $K$  in three space that 
we take to be an embedding (knot) or a curve with transversal self-intersections.  For this 
discussion, the Wilson loop will be denoted by the notation  $W_{K}(A) = <K|A>$ to 
denote the dependence on the loop $K$ and the field $A$.   It is usually indicated by the 
symbolism   $tr(Pe^{\oint_{K} A})$ .   Thus   
$$W_{K}(A) = <K|A>  = tr(Pe^{\oint_{K} A}).$$    Here the $P$  denotes  path ordered 
integration - we are integrating and exponentiating matrix valued functions, and so must 
keep track of the order of the operations.  The  symbol  $tr$  denotes the trace of the 
resulting matrix.
\vspace{3mm}

With the help of the Wilson loop functional on knots and links,  Witten  writes down a 
functional integral for link invariants in a 3-manifold  $M$:

$$Z(M,K) = \int DAe^{(ik/4 \pi)S(M,A)} tr(Pe^{\oint_{K} A}) $$

$$= \int DAe^{(ik/4 \pi)S}<K|A>.$$

\noindent
Here $S(M,A)$ is the Chern-Simons Lagrangian, as in the previous discussion. We 
abbreviate  $S(M,A)$  as $S$ and write  $<K|A>$  for the Wilson loop. Unless otherwise 
mentioned, the manifold  $M$  will be the three-dimensional sphere  $S^{3}$
\vspace{3mm} 

An analysis of the formalism of this functional integral  
reveals quite a bit about its role in knot theory.   This analysis depends upon key facts 
relating the curvature of the gauge field to both the Wilson loop and the Chern-Simons 
Lagrangian. The idea for using the curvature in this way is due to Lee Smolin 
\cite{Smolin} (See also \cite{Ramusino}). 
To this end, let us recall the local coordinate structure of the gauge field  $A(x)$,  where  
$x$  is a point in three-space.   We can 
write   $A(x)  =  A^{a}_{k}(x)T_{a}dx^{k}$  where  the index  $a$ ranges from $1$ to 
$m$ with the Lie
algebra basis $\{T_{1}, T_{2}, T_{3}, ..., T_{m}\}$.  The index $k$   goes from $1$  to  
$3$.     For each choice of $a$  and  $k$,  $A^{a}_{k}(x)$   is a smooth function defined 
on three-space.  
In  $A(x)$  we sum over the values of repeated indices.  The Lie algebra generators 
$T_{a}$  are  matrices  corresponding to a given representation of the Lie algebra of the 
gauge group $G.$   We assume some properties of these matrices as follows:
\vspace{3mm}

\noindent 1.  $[T_{a} , T_{b}] = i f^{abc}T_{c}$  where  $[x ,y] = xy - yx$ , and 
$f^{abc}$ 
(the matrix of structure constants)  is totally antisymmetric.  There is summation over 
repeated indices.
 \vspace{3mm}
 
\noindent 2.  $tr(T_{a}T_{b}) = \delta_{ab}/2$ where  $\delta_{ab}$ is the Kronecker 
delta  ($\delta_{ab} = 1$ if $a=b$ and zero otherwise).
\vspace{6mm}

We also assume some facts about curvature. (The reader may enjoy comparing with the 
exposition in \cite{K and P}.  But note the difference of conventions on the use of  $i$ in 
the Wilson loops and curvature definitions.)   The first fact  is the relation of Wilson loops 
and curvature for small loops:
 \vspace{3mm}

\noindent {\bf Fact 1.} The result of evaluating a Wilson loop about a very small planar 
circle around a point $x$ is proportional to the area enclosed by this circle times the 
corresponding value of the curvature tensor of the gauge field evaluated at $x$. The 
curvature tensor is  written  $$F^{a}_{rs}(x)T_{a}dx^{r}dy^{s}.$$  
 It is the local coordinate expression of  $F = dA +A \wedge A.$
\vspace{3mm} 

\noindent 
{\bf Application of Fact 1.}  Consider a given Wilson line  $<K|S>$.    
 Ask how its value will change if it is deformed infinitesimally in the neighborhood of a 
point $x$ on the line.  Approximate the change according to Fact 1, and regard the point 
$x$ as the place of curvature evaluation.  Let  $\delta<K|A>$  denote the change
 in the value of the line.    $\delta <K|A>$  is given by the formula    
 $$\delta <K|A> = dx^{r}dx^{s}F_{a}^{rs}(x)T_{a}<K|A>.$$  
 This  is the first order approximation to the change in the Wilson line.
 \vspace{3mm}

In this formula it  is understood that the Lie algebra matrices  $T_{a}$  are to be inserted 
into the Wilson line at the point $x$,  and that we are summing over repeated indices. This 
means that  each  $T_{a}<K|A>$ is  a new Wilson line obtained from the original  line  
$<K|A>$  by leaving the form of the loop unchanged,  but inserting the matrix  $T_{a}$ 
into that loop at the point  $x$.   In Figure 6 we have illustrated this mode of insertion of Lie algebra into the Wilson loop. Here and in further illustrations in this section we use $W_{K}(A)$ to denote the Wilson loop. Note that in the diagrammatic version shown in Figure 6 we have let small triangles with legs indicate $dx^{i}.$ The legs correspond to indices just as in our work in the last section with Lie algebras and chord diagrams. The curvature tensor is indicated as a circle with three legs corresponding to the indices of $F_{a}^{rs}.$
\vspace{3mm}  

\noindent
{\bf Notation.} In the diagrams in this section we have dropped mention of the factor of $(1/ 4 \pi)$ that occurs in the integral. This convention saves space in the figures. In these figures $L$ denotes the Chern--Simons Lagrangian. 
\vspace{3mm}  

\begin{figure}[htbp]
\vspace*{60mm}
\includegraphics{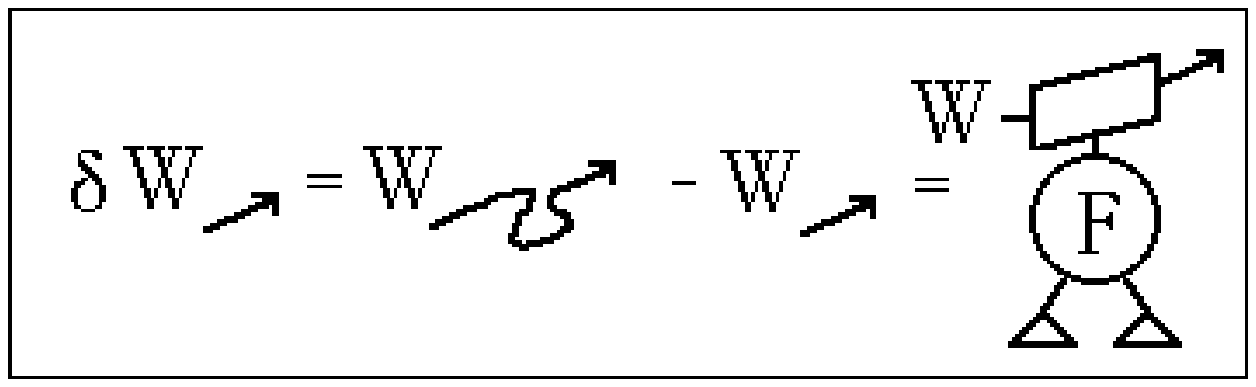}
\vspace*{13pt}
\begin{center}
{\bf Figure 6 --- Lie algebra and Curvature Tensor insertion into the Wilson Loop}
\end{center}
\end{figure}
\vspace{3mm}

\noindent {\bf Remark.}  In thinking about the Wilson line  
$<K|A> = tr(Pe^{\oint_{K} A})$,  it is helpful to recall Euler's formula for the 
exponential:  

$$e^{x} = lim_{n \rightarrow \infty}(1+x/n)^{n}.$$

\noindent
The  Wilson line  is  the limit, over partitions of the loop $K$,  of
products of the matrices  $(1 + A(x))$  where $x$ runs over the partition.  Thus we can 
write symbolically,

$$<K|A> =  \prod_{x \in K}(1 +A(x))$$  
$$=  \prod_{x \in K}(1 + A^{a}_{k}(x)T_{a}dx^{k}).$$

\noindent
It is understood that a product of matrices around a closed loop connotes the trace of the 
product.  The ordering is forced by the one dimensional nature of the loop.   Insertion of a 
given matrix into this product at a point on the loop is then a well-defined concept.   If  $T$  
is a given matrix then it is understood that   $T<K|A>$  denotes the insertion of $T$ into 
some point of the loop.  In the case above, it is understood from context in the formula that 
the insertion is to be performed at the point $x$  indicated in the argument of the curvature. 
 \vspace{3mm}
   
\noindent {\bf Remark.}  The   previous remark implies the following formula for the 
variation of the Wilson loop with respect to the gauge field:

$$\delta <K|A>/\delta (A^{a}_{k}(x))  =  dx^{k}T_{a}<K|A>.$$

\noindent
Varying the Wilson loop with respect to the gauge field results in the insertion of an 
infinitesimal Lie algebra element into the loop. Figure 7 gives a diagrammatic form for this formula. In that Figure we use a capital $D$ with up and down legs to denote the derivative $\delta /\delta (A^{a}_{k}(x)).$ Insertions in the Wilson line are indicated directly by matrix boxes placed in a representative bit of line.
\vspace{3mm}

\begin{figure}[htbp]
\vspace*{50mm}
\includegraphics{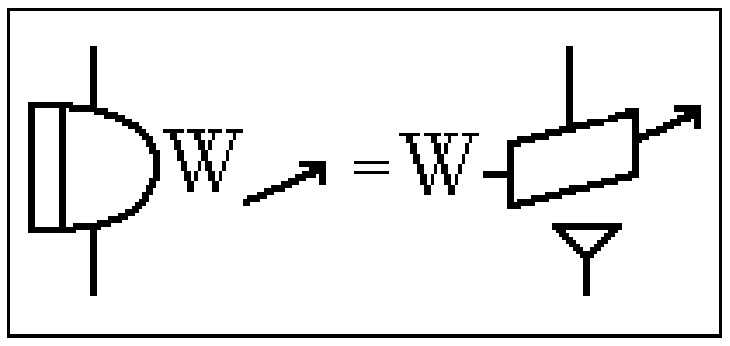}
\vspace*{13pt}
\begin{center}
{\bf Figure 7 --- Differentiating the Wilson Line}
\end{center}
\end{figure}
\vspace{3mm}

\noindent {\bf Proof.}  
$$\delta <K|A>/\delta (A^{a}_{k}(x))$$   

$$= \delta \prod_{y \in K}(1 + A^{a}_{k}(y)T_{a}dy^{k})/\delta (A^{a}_{k}(x))$$

$$= \prod_{y<x \in K}(1 + A^{a}_{k}(y)T_{a}dy^{k}) [T_{a}dx^{k}] \prod_{y>x \in 
K}(1 + A^{a}_{k}(y)T_{a}dy^{k})$$

$$= dx^{k}T_{a}<K|A>.$$ 
\vspace{3mm}

\noindent
{\bf Fact 2.}  The variation of the Chern-Simons Lagrangian  $S$  with respect to the 
gauge potential at a given point in three-space is related to the values of the curvature 
 tensor at that point by the following formula:

$$F^{a}_{rs}(x)  =  \epsilon_{rst} \delta S/\delta (A^{a}_{t}(x)).$$

\noindent
 Here  $\epsilon_{abc}$ is the epsilon symbol for three indices, i.e. it is $+1$ for positive 
permutations of $123$ and
$-1$ for negative permutations of $123$ and zero if any two indices are repeated.  A diagrammatic for this formula is shown in Figure 8.
\vspace{3mm}

\begin{figure}[htbp]
\vspace*{60mm}
\includegraphics{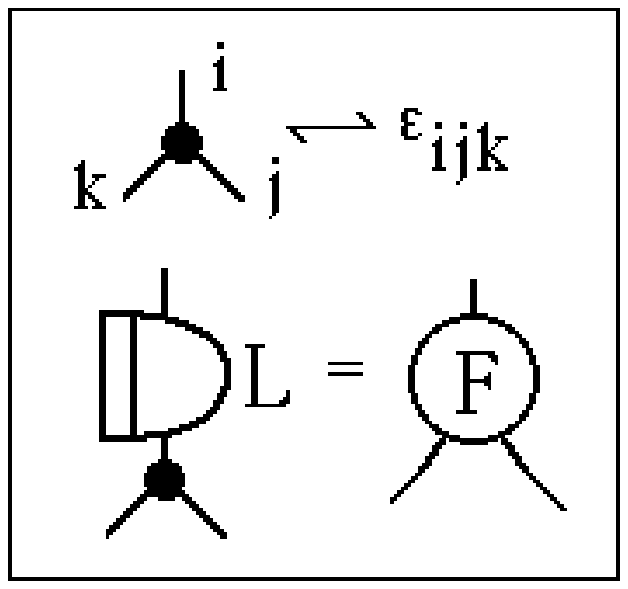}
\vspace*{13pt}
\begin{center}
{\bf Figure 8 --- Variational Formula for Curvature}
\end{center}
\end{figure}
\vspace{3mm}

With these facts at hand we are prepared to determine how the Witten integral behaves 
under a small deformation of the loop $K.$
 \vspace{3mm}

\noindent 
{\bf Theorem.}
1. Let   $Z(K) = Z(S^{3},K)$  and let  $\delta Z(K)$  denote the change of  
$Z(K)$ under an infinitesimal change in the loop  K.   Then

  $$ \delta Z(K) = (4 \pi i/k) \int dA e^{(ik/4\pi)S}[Vol] T_{a} T_{a} <K|A>$$

\noindent
 where $Vol = \epsilon_{rst} dx^{r} dx^{s} dx^{t}.$

The sum is taken over repeated indices, and the insertion is taken of the matrices   
$T_{a}T_{a}$  at the chosen point  $x$  on the loop $K$ that is regarded as the center 
of the deformation.  The volume element 
$Vol = \epsilon_{rst}dx_{r}dx_{s}dx_{t}$ is taken
with regard to the infinitesimal directions of the loop deformation from this point on the 
original loop.
\vspace{3mm}

\noindent 
2. The same formula applies, with a different interpretation,  to the case where  $x$  is  a 
double point of transversal self intersection of a loop K,  and the deformation consists in 
shifting one of the crossing segments perpendicularly to the plane of 
intersection  so that the self-intersection point disappears.  In this  case,  one  $T_{a}$  is 
inserted into each of the transversal crossing segments so that  $T_{a}T_{a}<K|A>$ 
denotes a Wilson loop with a self intersection  at  $x$   and insertions of $T_{a}$  at  $x + 
\epsilon_{1}$ and  $x + \epsilon_{2}$  where $\epsilon_{1}$ and $\epsilon_{2}$ denote 
small displacements along the two arcs of $K$ that intersect at $x.$  In this case, the 
volume form is nonzero, with two  directions coming from the plane of movement of one 
arc, and the perpendicular direction is the direction of the other arc.
 \vspace{3mm}

\noindent {\bf Proof.}
 
$$\delta Z(K)  =  \int DA e^{(ik/4 \pi)S} \delta <K|A>$$

$$= \int DA e^{(ik/4 \pi)S} dx^{r}dy^{s} F^{a}_{rs}(x) T_{a}<K|A>$$    

$$=  \int DA e^{(ik/4 \pi)S} dx^{r}dy^{s} \epsilon_{rst} (\delta S/\delta (A^{a}_{t}(x)))  T_{a}<K|A>$$         

$$= (-4 \pi i/k) \int DA  (\delta e^{(ik/4 \pi)S}/\delta (A^{a}_{t}(x))) \epsilon_{rst} 
dx^{r}dy^{s}T_{a}<K|A>$$

$$= (4 \pi i/k) \int DA  e^{(ik/4 \pi)S} \epsilon_{rst} dx^{r}dy^{s} (\delta 
T_{a}<K|A>/\delta (A^{a}_{t}(x)))$$

(integration by parts and the boundary terms vanish)

$$=  (4 \pi i/k) \int DA  e^{(ik/4 \pi)S}[Vol] T_{a}T_{a}<K|A>.$$

This completes the formalism of the proof.  In the case of part 2., a change of interpretation 
occurs at the point in the argument when the Wilson line is differentiated.  Differentiating a 
self intersecting Wilson line at a point of self intersection is equivalent to differentiating the 
corresponding product of
matrices with respect to a variable that occurs at two points in the product (corresponding 
to the two places where the loop passes through the point).  One of these derivatives gives 
rise to a term with volume form equal to zero, the other term is the one that is described in 
part 2.  This completes the proof of the Theorem.  //
\vspace{3mm}
 
\noindent
The formalism of this proof is illustrated in Figure 9.
\vspace{3mm}

\begin{figure}[htbp]
\vspace*{170mm}
\includegraphics{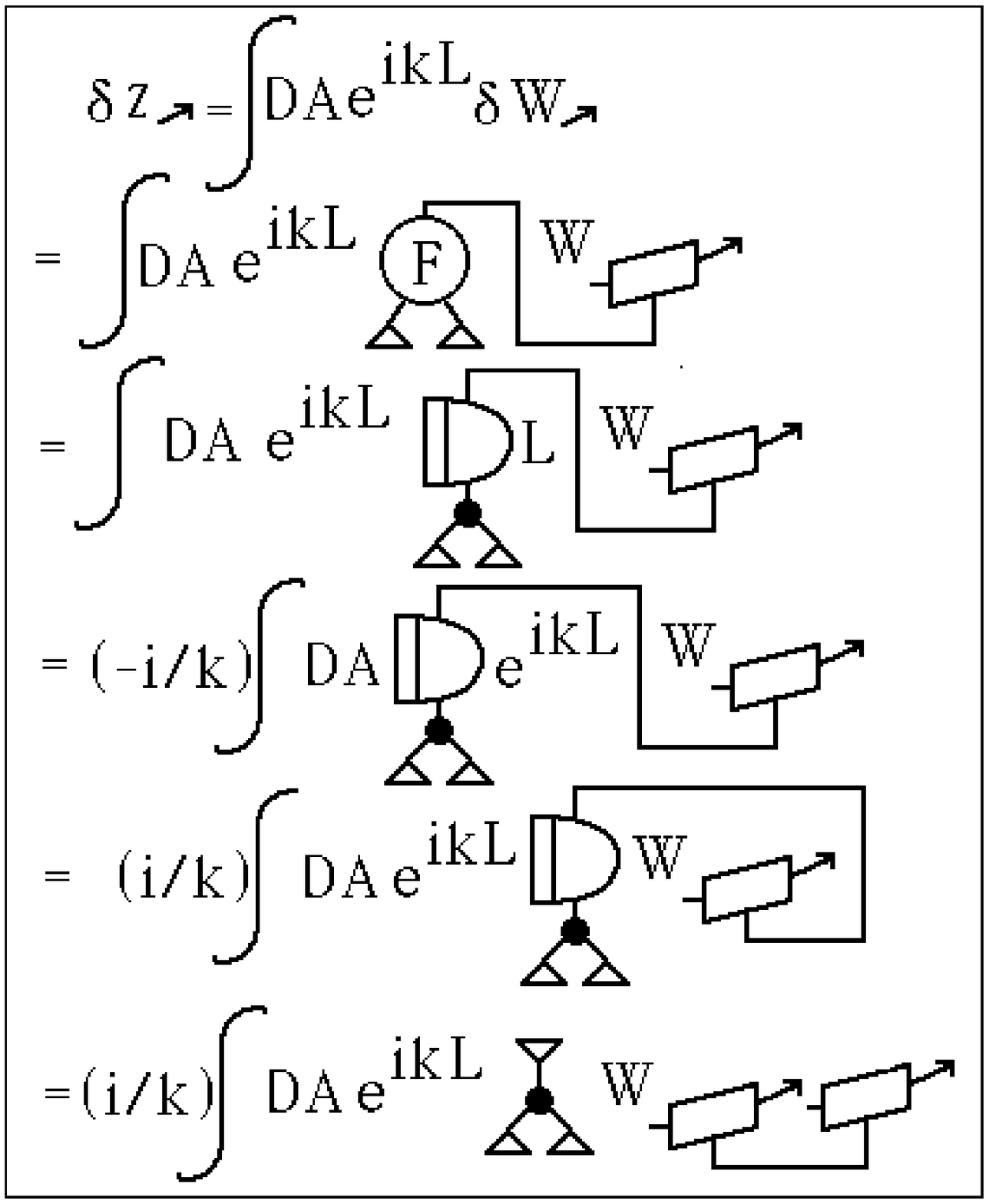}
\vspace*{13pt}
\begin{center}
{\bf Figure 9 --- Varying the Functional Integral by Varying the Line}
\end{center}
\end{figure}
\vspace{3mm}

In the case of switching a crossing the key point is to write the crossing switch as a 
composition of first moving a segment to obtain a transversal intersection of the diagram 
with itself, and then to continue the motion to complete the switch.  One then analyses 
separately the case where $x$  is  a double point of transversal self intersection of a loop 
$K,$  and the deformation consists in shifting one of the crossing segments 
perpendicularly to the plane of 
intersection  so that the self-intersection point disappears.  In this  case,  one  $T_{a}$  is 
inserted into each of the transversal crossing segments so that  $T^{a}T^{a}<K|A>$ 
denotes a Wilson loop with a self intersection  at  $x$   and insertions of $T^{a}$  at  
$x + \epsilon_{1}$ and  $x + \epsilon_{2}$   as in part $2.$ of the Theorem above. The 
first insertion is in the moving line, due to curvature. The second insertion is the 
consequence of differentiating the self-touching Wilson line. Since this line can be regarded 
as a product, the differentiation occurs twice at the point of intersection, and it is the second 
direction that produces the non-vanishing volume form.
 \vspace{3mm}

Up to the choice of our conventions for constants, the switching formula is, as shown 
below (See Figure 10).

 $$Z(K_{+}) -  Z(K_{-}) =  (4 \pi i/k)\int DA  e^{(ik/4\pi)S} T_{a}T_{a}<K_{**}|A>$$ 
$$= (4 \pi i/k) Z(T^{a}T^{a}K_{**}),$$

\noindent
where $K_{**}$ denotes the result of replacing the crossing by a self-touching crossing. 
We distinguish this from adding a graphical node at this crossing by using the double star 
notation.
\vspace{3mm}

\begin{figure}[htbp]
\vspace*{70mm}
\includegraphics{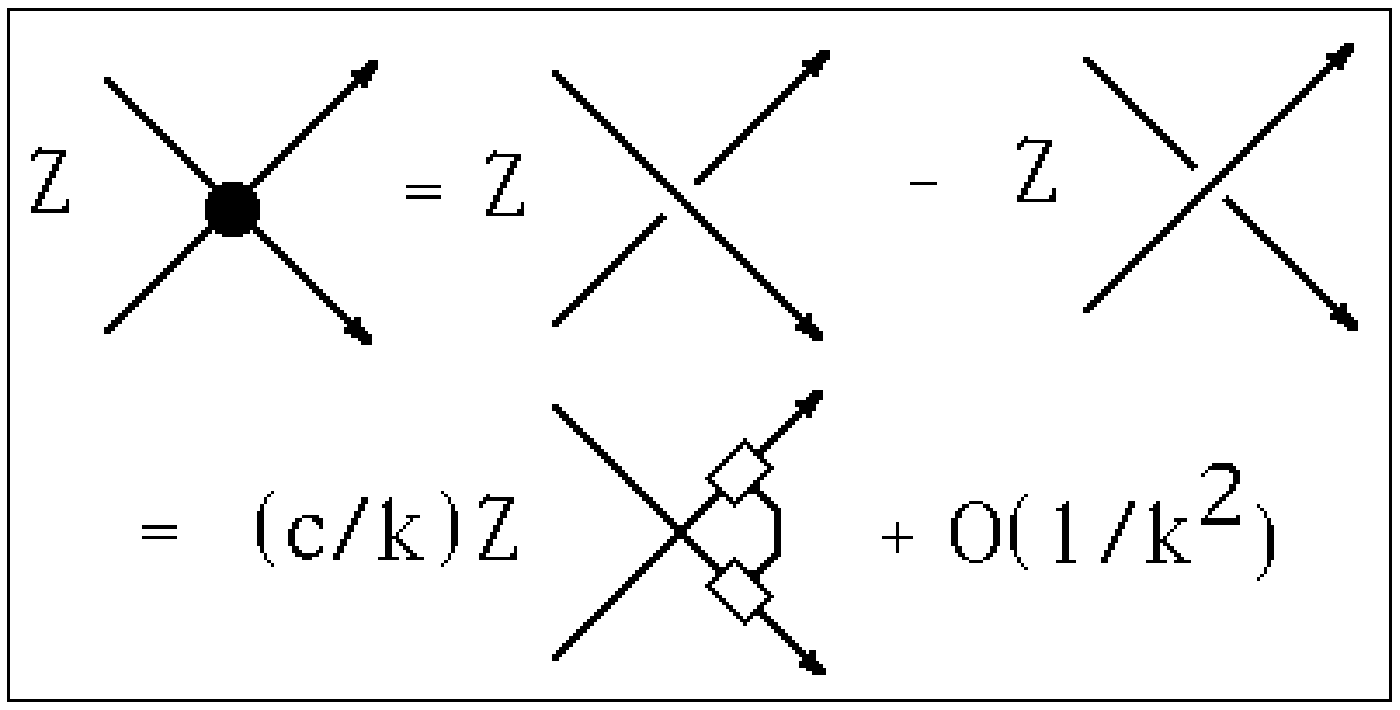}
\vspace*{13pt}
\begin{center}
{\bf Figure 10 --- The Difference Formula}
\end{center}
\end{figure}
\vspace{3mm}

A key point is to notice that the Lie algebra insertion for this difference is exactly what is 
done (in chord diagrams) to make the weight systems for Vassiliev invariants (without the 
framing compensation).  Here we take formally the perturbative expansion of the Witten 
integral to obtain Vassiliev invariants as coefficients of the powers of ($1/k^{n}$). Thus 
the formalism of the Witten functional integral takes one directly to these weight systems in 
the case of the classical Lie algebras. In this way the functional integral  is central to the 
structure of the Vassiliev invariants. 
\vspace{3mm}

\subsection{The Loop Transform}

Suppose that $\psi (A)$ is a (complex valued) function defined on gauge fields. Then we define formally the {\em loop transform}
$\widehat{\psi}(K)$, a function on embedded loops in three dimensional space, by the formula

$$\widehat{\psi}(K) = \int DA \psi(A) W_{K}(A).$$

\noindent
If $\Delta$ is a differential operator defined on $\psi(A),$ then we can use this integral transform to shift the effect of $\Delta$ to an operator on loops via integration by parts:

$$\widehat{ \Delta \psi }(K) = \int DA \Delta \psi(A) W_{K}(A)$$

$$ = - \int DA  \psi(A) \Delta W_{K}(A).$$

\noindent
When $\Delta$ is applied to the Wilson loop the result can be an understandable geometric or topological operation.  In Figures 11, 12 and 13 we illustrate this situation with diagrammatically defined operators $G$ and $H.$ 
\vspace{3mm}

\begin{figure}[htbp]
\vspace*{160mm}
\includegraphics{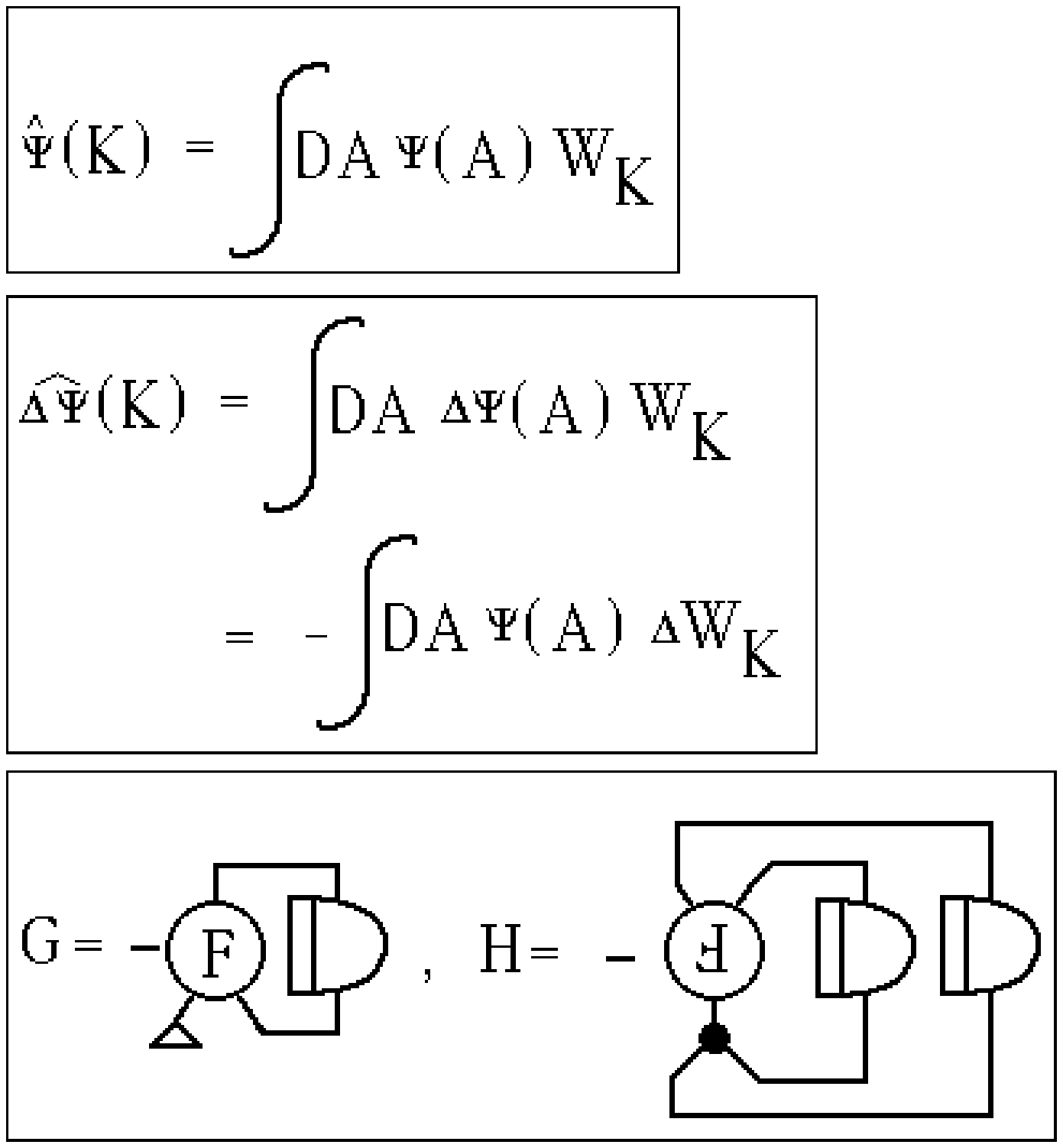}
\vspace*{13pt}
\begin{center}
{\bf Figure 11--- The Loop Transform and Operators G and H}
\end{center}
\end{figure}
\vspace{3mm}

\begin{figure}[htbp]
\vspace*{160mm}
\includegraphics{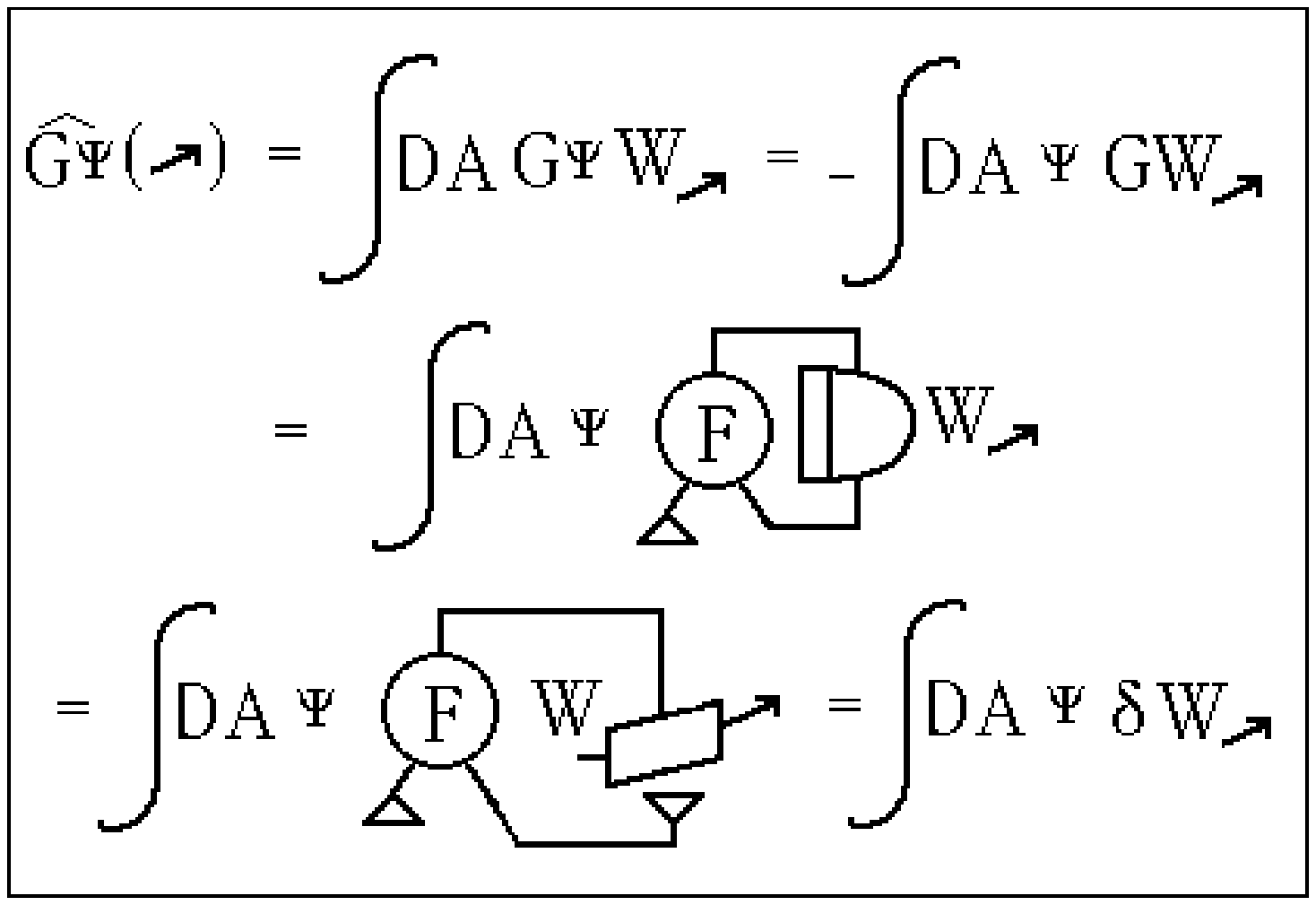}
\vspace*{13pt}
\begin{center}
{\bf Figure 12 --- The Diffeomorphism Constraint}
\end{center}
\end{figure}
\vspace{3mm}

\begin{figure}[htbp]
\vspace*{160mm}
\includegraphics{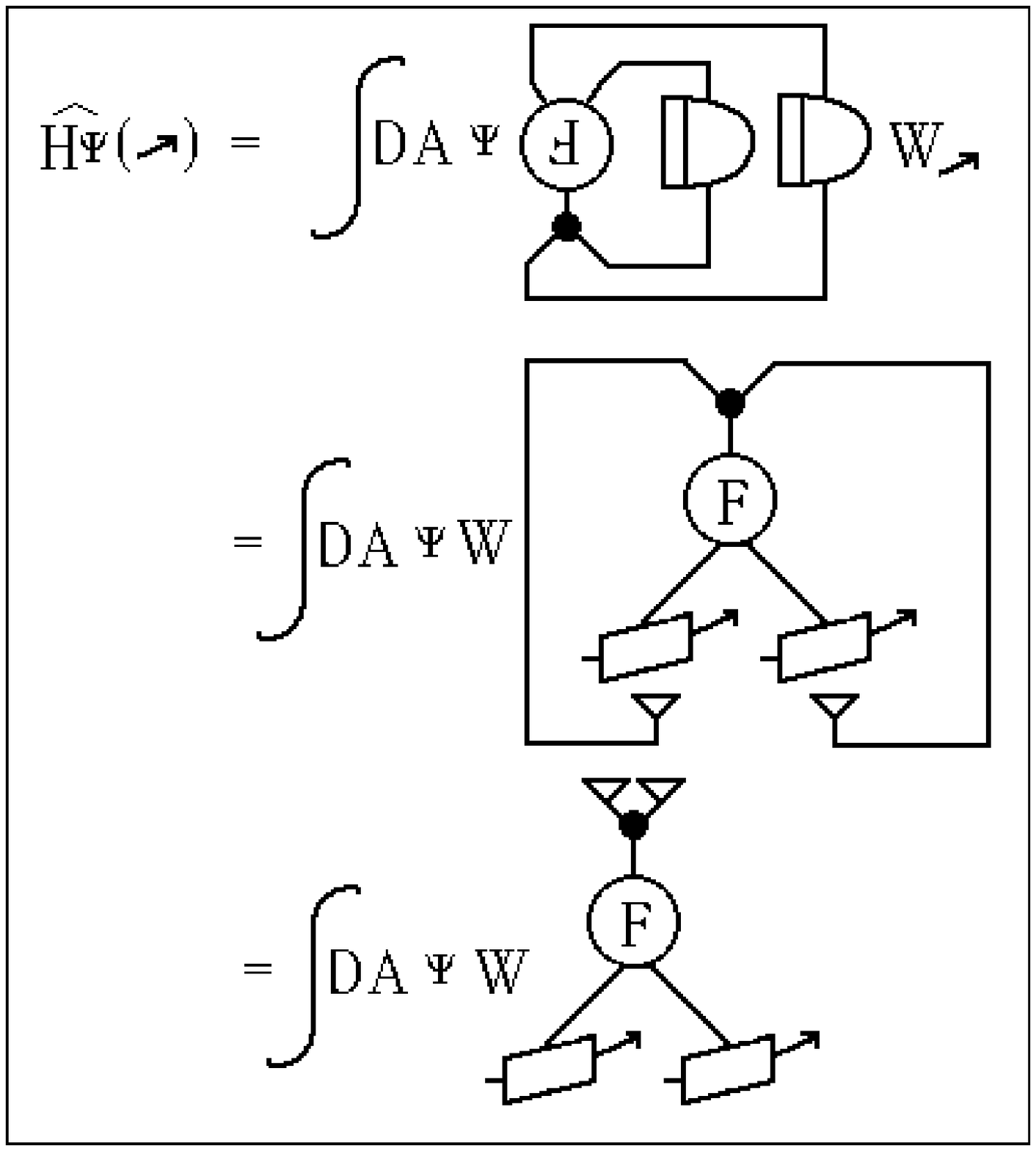}
\vspace*{13pt}
\begin{center}
{\bf Figure 13 --- The Hamiltonian Constraint}
\end{center}
\end{figure}
\vspace{3mm}

\noindent
We see from Figure 12 that 

$$\widehat{ G \psi }(K) =  \delta \widehat{ \psi }(K)$$

\noindent
where this variation refers to the effect of  varying $K$ by a small loop. As we saw in this section, this means that if  $\widehat{ \psi }(K)$  is a topological invariant of knots and links, then 
$\widehat{ G \psi }(K) =0$ for all embedded loops $K.$  This condition is a transform analogue of the equation $G \psi(A) =0.$
This equation is the differential analogue of an invariant of knots and links. It may happen that $\delta \widehat{ \psi }(K)$ is not strictly zero, as in the case of our framed knot invariants.
For example with $$\psi(A) =  e^{(ik/4\pi) \int tr(A \wedge dA + (2/3)A \wedge A \wedge A)}$$
we conclude that $\widehat{ G \psi }(K)$ is zero for flat deformations (in the sense of this section) of the loop $K,$ but can be non-zero in the presence of a twist or curl. In this sense the loop transform provides a subtle variation on the strict condition $G \psi(A) =0.$
\vspace{3mm}

In \cite{ASR} and earlier publications by these authors, the loop transform is used to study a reformulation and quantization of Einstein gravity.  The differential geometric gravity theory is reformulated in terms of a background gauge connection and in the quantization, the Hilbert space consists in functions $\psi(A)$ that are required to satisfy the constraints 

$$G \psi =0$$
 
\noindent
and 

$$H \psi =0$$

\noindent
where $H$ is the operator shown in Figure 13. Thus we see that 
$\widehat{G}(K)$ can be partially zero in the sense of producing a framed knot invariant, and (from Figure 13 and the antisymmetry of the epsilon) that $\widehat{H}(K)$ is zero for non-self intersecting loops.  This means that the loop transforms of $G$ and $H$ can be used to investigate a subtle variation of the original scheme for the quantization of gravity. This program is being actively pursued by a number of researchers. The Vassiliev invariants arising from a topologically invariant loop transform should be of significance to this theory.  This theme will be explored in a subsequent paper. 
\vspace{3mm}

\section{Wilson Lines, Axial Gauge and the Kontsevich Integrals}

In this section we follow the gauge fixing method used by Fr\"ohlich and King 
\cite{Frohlich and King}.  Their paper was written before the advent of Vassiliev 
invariants, but contains, as we shall see, nearly the whole story about the Kontsevich 
integral.  A similar approach to ours can be found in \cite{LP}. In our case we have simplified the determination of the inverse operator for this formalism and we have given a few more details about the calculation of the correlation functions than is customary in physics literature.  I hope that this approach makes this subject more accessible to mathematicians. A heuristic argument of this kind contains a great deal of valuable mathematics. It is clear that these matters will eventually be given a fully rigorous treatment. In fact, in the present case there is a rigorous treatment, due to Albevario and Sen-Gupta \cite{AS} of the functional integral {\em after} the light-cone  gauge has been imposed.
\vspace{3mm}

\noindent
Let $(x^{0}, x^{1}, x^{2})$ denote a point in three dimensional space.
Change to light-cone coordinates
$$x^{+} = x^{1} + x^{2}$$ and
$$x^{-} = x^{1} - x^{2}.$$

\noindent
Let $t$ denote $x^{0}.$
\vspace{3mm}

\noindent
Then the gauge connection can be written in the form
$$A(x) = A_{+}(x)dx^{+} + A_{-}(x)dx^{-} + A_{0}(x)dt.$$
\vspace{3mm}

\noindent
Let $CS(A)$ denote the Chern-Simons integral (over the three dimensional sphere)

$$CS(A) = (1/4\pi)\int tr(A \wedge dA + (2/3) A \wedge A \wedge A).$$

\noindent
We define {\em axial gauge} to be the condition that $A_{-} = 0.$
We shall now work with the functional integral of the previous section under the axial 
gauge restriction. In axial gauge we have that 
$A \wedge A \wedge A = 0$ and so 

$$CS(A) = (1/4\pi)\int tr(A \wedge dA).$$

\noindent
Letting $\partial_{\pm}$ denote partial differentiation with respect to $x^{\pm}$, we get 
the following formula in axial gauge
$$A \wedge dA = (A_{+} \partial_{-} A_{0} - A_{0} \partial_{-}A_{+})dx^{+} \wedge 
dx^{-} \wedge dt.$$

\noindent
Thus, after integration by parts, we obtain the following formula for the Chern-Simons 
integral:

 $$CS(A) = (1/2 \pi) \int tr(A_{+} \partial_{-} A_{0}) dx^{+} \wedge dx^{-} \wedge 
dt.$$

\noindent
Letting $\partial_{i}$ denote the partial derivative with respect to $x_{i}$, we have that 
$$\partial_{+} \partial_{-} = \partial_{1}^{2} - \partial_{2}^{2}.$$  If we replace 
$x^{2}$ with $ix^{2}$ where 
$i^{2} = -1$, then $\partial_{+} \partial_{-}$  is replaced by
$$\partial_{1}^{2} + \partial_{2}^{2} = \nabla^{2}.$$ 
We now make this replacement so that the analysis can be expressed over the complex 
numbers.
\vspace{3mm}

\noindent
Letting $$z = x^{1} + ix^{2},$$ it is well known that 
$$\nabla^{2} ln(z) = 2 \pi \delta(z)$$
where $\delta(z)$ denotes the Dirac delta function and $ln(z)$ is the natural logarithm of 
$z.$  Thus we can write
$$(\partial_{+} \partial_{-})^{-1} = (1/2 \pi)ln(z).$$
Note that $\partial_{+} = \partial_{z} = \partial /\partial z$  after the replacement of 
$x^{2}$ by $ix^{2}.$   As a result we have that 
$$(\partial_{-})^{-1} = \partial_{+} (\partial_{+} \partial_{-})^{-1} =
\partial_{+} (1/2 \pi)ln(z) = 1/2 \pi z.$$

\noindent
Now that we know the inverse of the operator $\partial_{-}$ we are in a position to treat 
the Chern-Simons integral as a quadratic form in the pattern 
$$ (-1/2)<A, LA> = - iCS(A)$$  
where the operator 
$$L =  \partial_{-}.$$
Since we know $L^{-1}$, we can express the functional integral as a Gaussian integral:
\vspace{3mm}

\noindent
We replace 

$$Z(K) = \int DAe^{ikCS(A)} tr(Pe^{\oint_{K} A})$$ 

\noindent
by

$$Z(K) = \int DAe^{iCS(A)} tr(Pe^{\oint_{K} A/\sqrt k})$$ 

\noindent
by sending $A$ to $(1/ \sqrt k)A$. We then replace this version by

$$Z(K) = \int DAe^{(-1/2)<A, LA>} tr(Pe^{\oint_{K} A/\sqrt k}).$$

\noindent
In this last formulation we can use our knowledge of $L^{-1}$ to determine the the 
correlation functions and express $Z(K)$ perturbatively in powers of $(1/ \sqrt k).$
\vspace{3mm}

\noindent
{\bf Proposition.}
Letting 
$$<\phi(A)> = \int DA e^{(-1/2)<A, LA>}\phi(A) / \int DA e^{(-1/2)<A, LA>}$$  for 
any functional $\phi(A)$,
we find that
$$<A_{+}^{a}(z,t)A_{+}^{b}(w,s)> = 0,$$
$$<A_{0}^{a}(z,t)A_{0}^{b}(w,s)> = 0,$$
$$<A_{+}^{a}(z,t)A_{0}^{b}(w,s)> = \kappa \delta^{ab} \delta(t-s)/(z-w)$$  where 
$\kappa$ is a constant.
\vspace{3mm}

\noindent
{\bf Proof Sketch.}
Let's recall how these correlation functions are obtained.
The basic formalism for the Gaussian integration is in the pattern

$$<A(z)A(w)> = \int DA e^{(-1/2)<A, LA>} A(z)A(w) / \int DA e^{(-1/2)<A, LA>}$$

$$ =  ((\partial / \partial J(z)) (\partial / \partial J(w)) |_{J=0})
 e^{(1/2)<J, L^{-1}J>}$$

\noindent
Letting $G*J(z) = \int dw G(z-w)J(w)$, we have that when

$$LG(z) = \delta(z)$$  

\noindent
($\delta(z)$ is a Dirac delta function of $z$.) then

$$LG*J(z) = \int dw LG(z-w)J(w) = \int dw \delta(z-w) J(w) = J(z)$$

\noindent
Thus $G*J(z)$ can be identified with $L^{-1}J(z)$.
\vspace{3mm}

\noindent
In our case 
$$G(z) = 1/ 2 \pi z$$  

\noindent
and 
 
$$L^{-1}J(z) = G*J(z) = \int dw J(w)/(z-w).$$

\noindent
Thus 

$$<J(z),L^{-1}J(z)> = <J(z), G*J(z)> = (1/ 2\pi) \int tr(J(z) (\int dw J(w)/(z-w)) dz$$

$$ = (1/ 2\pi) \int \int dz dw tr(J(z)J(w))/(z-w).$$

\noindent
The results on the correlation functions then follow directly from differentiating this 
expression.  Note that the Kronecker delta on Lie algebra indices is a result of the corresponding Kronecker delta in the trace formula
$tr(T_{a}T_{b}) = \delta_{ab}/2$ for products of Lie algebra generators. The Kronecker delta for the $x^{0} = t, s$ coordinates is a consequence of the evaluation at $J$ equal to zero.//
\vspace{3mm}

We are now prepared to give an explicit form to the perturbative expansion for 

$$<K>= Z(K)/\int DAe^{(-1/2)<A, LA>}$$

$$= \int DAe^{(-1/2)<A, LA>} tr(Pe^{\oint_{K} A/\sqrt k})/ \int DAe^{(-1/2)<A, LA>}$$

$$ = \int DAe^{(-1/2)<A, LA>} tr(\prod_{x \in K} (1 +  (A/\sqrt k)))/\int DAe^{(-
1/2)<A, LA>}$$

$$ = \sum_{n} (1/k^{n/2}) \oint_{K_{1} < ... < K_{n}} <A(x_{1}) ... A(x_{n})>.$$

\noindent
The latter summation can be rewritten (Wick expansion) into a sum over products of pair 
correlations, and we have already worked out the values of these. In the formula above we 
have written  $K_{1} < ... < K_{n}$ to denote the integration over variables $x_{1} , ... 
x_{n}$ on $K$ so that $x_{1} <  ... < x_{n}$ in the ordering induced on the loop $K$ by 
choosing a basepoint on the loop.  After the Wick expansion, we get

$$<K> = \sum_{m} (1/k^{m}) \oint_{K_{1} < ... < K_{n}} 
\sum_{P= \{x_{i} < x'_{i}| i = 1, ... m\}}
\prod_{i}<A(x_{i})A(x'_{i})>.$$

\noindent
Now we know that 

$$<A(x_{i})A(x'_{i})> = 
<A^{a}_{k}(x_{i})A^{b}_{l}(x'_{i})>T_{a}T_{b}dx^{k}dx^{l}.$$

\noindent
Rewriting this in the complexified axial gauge coordinates, the only contribution is 

$$<A_{+}^{a}(z,t)A_{0}^{b}(s,w)> = \kappa \delta^{ab} \delta(t-s)/(z-w).$$

\noindent
Thus

$$<A(x_{i})A(x'_{i})>$$ 
$$=<A^{a}_{+}(x_{i})A^{a}_{0}(x'_{i})>T_{a}T_{a}dx^{+} \wedge dt + 
<A^{a}_{0}(x_{i})A^{a}_{+}(x'_{i})>T_{a}T_{a}dx^{+} \wedge dt$$
$$= (dz-dz')/(z-z') [i/i']$$
where $[i/i']$ denotes the insertion of the Lie algebra elements
$T_{a}T_{a}$ into the Wilson loop.
\vspace{3mm}

\noindent
As a result, for each partition of the loop and choice of pairings 
$P= \{x_{i} < x'_{i}| i = 1, ... m\}$ we get an evaluation $D_{P}$ of the trace of these 
insertions into the loop. This is the value of the corresponding chord diagram in the weight 
systems for Vassiliev invariants. These chord diagram evaluations then figure in our 
formula as shown below:

$$<K> = \sum_{m} (1/k^{m}) \sum_{P} D_{P} \oint_{K_{1} < ... < K_{n}}
\bigwedge_{i=1}^{m}(dz_{i} - dz'_{i})/((z_{i} - z'_{i})$$

\noindent
This is a Wilson loop ordering version of the Kontsevich integral. To see the usual form of 
the integral appear, we change from the time variable (parametrization) associated with the 
loop itself to  time variables associated with a specific global direction of time in three 
dimensional space that is perpendicular to the complex plane defined by the axial gauge 
coordinates. It is easy to see that this results in one change of sign for each segment of the 
knot diagram supporting a pair correlation where the segment is oriented (Wilson loop 
parameter)  downward with respect to the global time direction. This results in the rewrite 
of our formula to

$$<K> = \sum_{m} (1/k^{m}) \sum_{P} (-1)^{|P \downarrow |} D_{P} \int_{t_{1} < 
... < t_{n}}
\bigwedge_{i=1}^{m}(dz_{i} - dz'_{i})/((z_{i} - z'_{i})$$

\noindent
where $|P \downarrow |$ denotes the number of points $(z_{i},t_{i})$ or $(z'_{i},t_{i})$  
in the pairings where the knot diagram is  oriented downward with respect to global time. 
The integration around the Wilson loop has been replaced by integration in the vertical time 
direction and is so indicated by the replacement of $\{ K_{1} < ... < K_{n} \}$ with $\{ 
t_{1} < ... < t_{n} \}$
\vspace{3mm}

\noindent
The coefficients of $1/k^{m}$  in this expansion are exactly the Kontsevich integrals for 
the weight systems $D_{P}$. See Figure 14.
\vspace{3mm}

\begin{figure}[htbp]
\vspace*{160mm}
\includegraphics{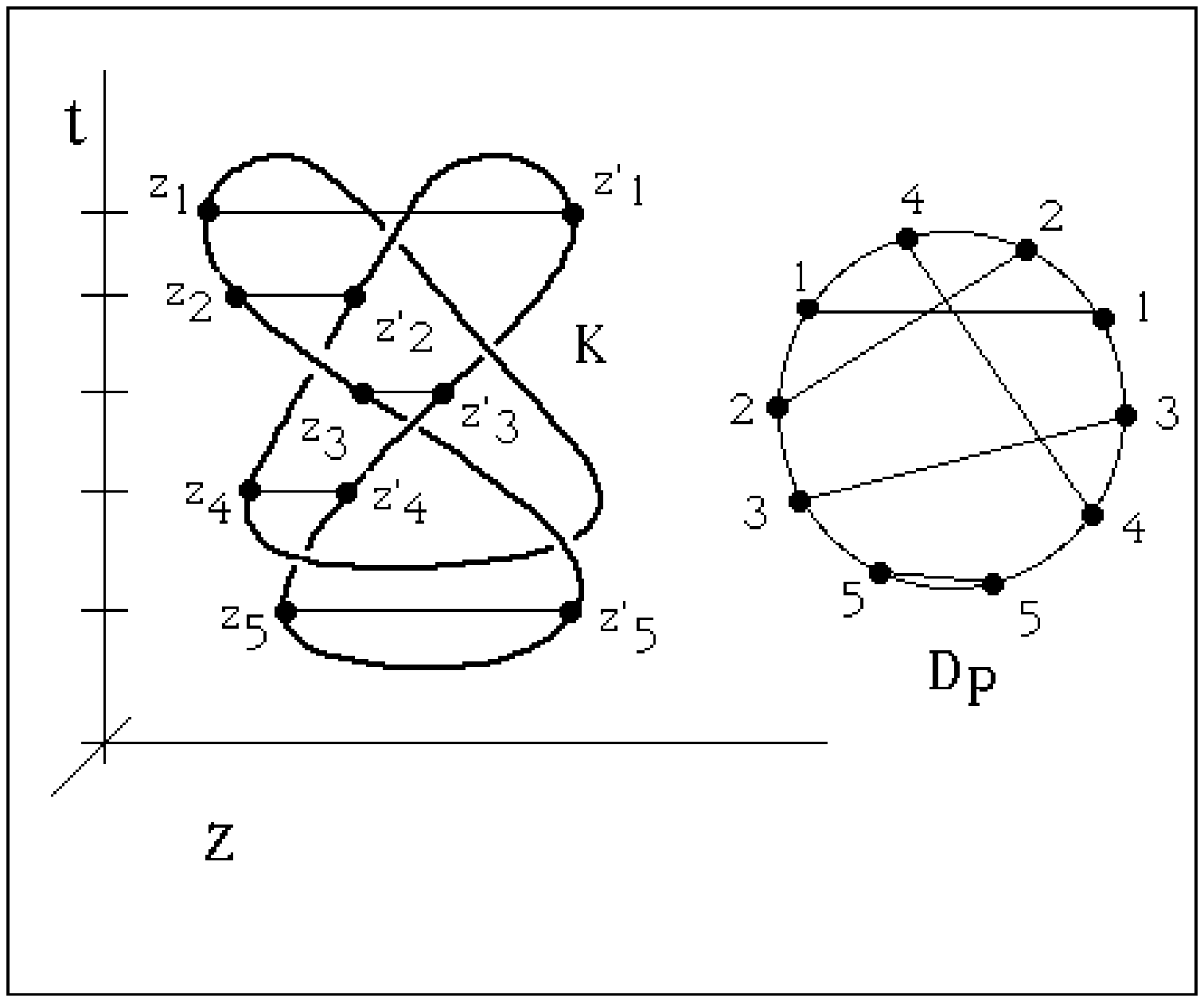}
\vspace*{13pt}
\begin{center}
{\bf Figure 14 --- Applying The Kontsevich Integral}
\end{center}
\end{figure}
\vspace{3mm}

\noindent
It was Kontsevich's insight to see (by different means) that 
these integrals could be used to construct Vassiliev invariants from arbitrary weight 
systems satisfying the four-term relations.  Here we have seen how these integrals arise 
naturally in the axial gauge fixing of the Witten functional integral. 
\vspace{3mm}

\noindent
{\bf Remark.}  The careful reader will note that we have not made a discussion of the role of the maxima and minima of the space curve of the knot with respect to the height direction ($t$).  In fact one has to take these maxima and minima very carefully into account and to divide by the corresponding evaluated loop pattern (with these maxima and minima) to make the Kontsevich integral well-defined and actually invariant under ambient  isotopy (with appropriate framing correction as well). The corresponding difficulty appears here in the fact that because of the gauge choice the Wilson lines are actually only defined in the complement of the maxima and minima and one needs to analyse a limiting procedure to take care of the inclusion of these points in the Wilson line. This points to one of the places where this correspondence with the Kontsevich integrals as Feynman integrals for Witten's functional integral could stand closer mathematical scrutiny. One purpose of this paper has been to outline the correspondences that exist and to put enough light on the situation to allow a full story to eventually appear.

 \end{document}